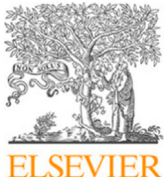

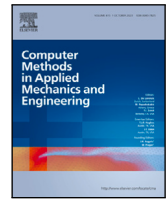

# A variationally consistent membrane wrinkling model based on spectral decomposition of the strain tensor

Daobo Zhang, Josef Kiendl *

*Institute for Engineering Mechanics and Structural Analysis - University of the Bundeswehr Munich, Werner-Heisenberg-Weg 39, 85579, Neubiberg, Germany*

## ARTICLE INFO

## ABSTRACT

We propose a novel variationally consistent membrane wrinkling model for analyzing the mechanical responses of wrinkled thin membranes. The elastic strain energy density is split into tensile and compressive terms via a spectral decomposition of the strain tensor. Tensile and compressive parts of the stress and constitutive tensors are then obtained via consistent variation from the respective strain energies. Considering only the positive part of the strain energy in the variational formulation, we obtain a membrane with zero compressive stiffness. By adding the negative strain energy multiplied with a very small factor, we further obtain a residual compressive stiffness, which improves stability and allows handling also states of slackening. Comparison with results from analytical, numerical and experimental examples from the literature on membrane wrinkling problems demonstrate the great performance and capability of the proposed approach, which is also compatible with commercial finite element software.

## 1. Introduction

Wrinkling is a common phenomenon in elastic membrane structures, which can significantly affect their mechanical properties and performance. Related applications include flexible floating structures, parachutes, solar sails, airbags, as well as human skin and tissues, among others. The investigation of wrinkling in thin membranes with diverse methods has been ongoing for several decades to enhance the understanding of their underlying characteristics. Despite extensive research, modeling and predicting the behavior of wrinkled membranes remains challenging. This is primarily because membranes exhibit minimal resistance to compression and bending. As a result, most compressive stresses within membranes are released through localized instability. It leads to out-of-plane displacements, which are typically referred to as wrinkles [1]. In this context, three possible states of a membrane are typically considered, see Fig. 1. A membrane is taut, when it is subjected to tension in two orthogonal directions, wrinkled, when it is subjected to tension in one direction and compression in the other one, and slack, when it is subjected to compression in both directions. There exist different approaches to determine the state of a membrane, namely the principal stress criterion, the principal strain criterion, and the mixed criterion. In the first one, principal stresses are used to distinguish between tension and compression, in the second one this is done based on the principal strains. The mixed criterion is a combination of these two, combining the advantages of both. It has been discussed by several authors that both the stress-based and strain-based criteria can lead to incorrect judgment of the membrane state in specific cases, and that the mixed criterion can be considered as the physically most correct one, see, e.g., [2–4]. Nevertheless, also the stress- and strain-based criteria can be considered as reasonable criteria and are used in as basis for different wrinkling models in the literature, e.g., [5–8].

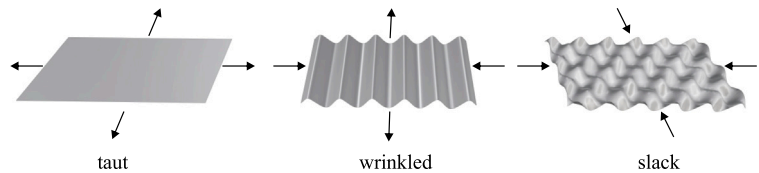

taut                                    wrinkled                                    slack

**Fig. 1.** Possible states of a membrane: taut, wrinkled, slack.

Several analytical solutions have been developed to address membrane wrinkling issues, as evidenced by the works in [9–12]. These solutions are derived within the tension field theory (TFT) framework, a theoretical approach that describes the highly buckled or wrinkled state of membranes undergo certain in-plane displacements that exceed the critical thresholds necessary for initiating buckling phenomena. Notably, Stein and Hedgepeth [5] derived a theory at small strains to predict the stresses and average deformations of partly wrinkled membranes by introducing an unknown function termed variable Poisson's ratio $\lambda$ and vanishing one of the principal stresses in wrinkled regions. Additionally, they provided analytical solutions for various partly wrinkled membrane problems, commonly employed as benchmarks to validate the stresses predicted by newly developed methods. Other analytical approaches are mentioned for completeness, and those are based on a geometrically nonlinear Föppl–von Kármán theory [13–15] and can depict the wrinkled configuration. However, these analytical solutions are often limited to specific problems and may lack applicability for cases that cannot be solved analytically. Therefore, numerical methods are preferred for analyzing wrinkled membrane problems, owing to their flexibility and ability to handle complex geometries and boundary conditions.

Previous numerical methods treating wrinkling in the finite element analysis framework can be broadly classified into two primary categories. This classification is determined by the requirement to obtain wrinkle details, including the number of wrinkles, wavelength, amplitude, wrinkle pattern, etc. The first relies on the buckling theory, which involves utilizing a highly dense mesh of shell elements with bending stiffness to resolve wrinkle details explicitly [16,17]. However, this method can be computationally expensive due to the need for a sufficiently refined mesh, and as pointed out in [18], the wrinkle length scale is frequently inconsistent with that of a finite element, leading to a mesh-dependent solution. In addition, a geometric imperfection must often be imposed on the membrane to induce wrinkles, which may also affect the results [19]. Consequently, the costs of this approach can rapidly become prohibitive, thereby limiting its practicality. Other relevant studies refer to Verhelst et al. [20], Flores and Oñate [21] and Fu et al. [22], to name a few.

The second one is to use a coarse mesh of membrane elements embedded with appropriate wrinkling models, mainly founded upon the assumptions of tension field theory [23]. It assumes that the wrinkled membrane exhibits zero compressive stiffness. Following the wrinkling process, the stress field within the wrinkled regions is presumed to be in a uniaxial tension stress state, whereby only one of the principal stress components is non-zero. Although the detailed wrinkles cannot be captured with these assumptions, it enables the identification of wrinkled regions and the directions of wrinkles, providing valuable insights into the overall behavior of the wrinkled membrane. Our work presented in this paper is carried out by following the second strategy. Thus, a literature review of wrinkling models is given below.

Most tension field theory-based wrinkling models can be broadly classified into two groups: kinematic modification methods and material modification methods. In general, the concept of the kinematic method is to modify the deformation gradient tensor by introducing additional parameters. Remarkably, the Roddeman model [24,25] has garnered significant attention, which extends the earliest Wu model [26] to anisotropic materials by incorporating a virtual elongation in the deformation tensor to satisfy the uniaxial tension condition. Among the subsequent developments of this model [27–31], Nakashino and Natori [29] circumvented the cumbersome linearization of the modified strain tensor to adjust the stress–strain tensor at the element level and recently Nakshino's model has been applied in isogeometric analysis [32]. Other models, such as Myazaki's model [1], which incorporates virtual elongation and shear deformation simultaneously into the deformation tensor, and direct modification of Green–Lagrange strain tensor [33–35], can also be classified as kinematic modification methods.

Material modification methods, which modify the stress–strain relation, can be divided into four subcategories. The first class is to modify the material parameters. For instance, Miller and Hedgepeth [7] and Miller et al. [8] extended the variable Poisson's ratio concept [5] into finite element static analysis and proposed the iterative material properties (IMP) model, where the local elasticity matrix is iteratively revised according to the state of strain in a previous load increment. The second and third ones are to use penalization techniques [2,3,6,18,36], and projection matrices [37–39] to eliminate compressive stresses by softening the constitutive tensor. However, the methods mentioned earlier require an explicit wrinkling criterion to evaluate the membrane states such as taut, wrinkled, and slack, as shown in Fig. 1. This reliance may lead to convergence problems due to the sudden changes in the tangent stiffness matrix. Thus, these shortcomings have motivated the fourth approach based on Pipkin's method, as documented in [40–44]. In this method, the strain energy density is replaced by a relaxed energy density acquired through an optimization procedure. Consequently, the stresses derived from the relaxed energy density are comparable to those obtained via the tension field theory. However, the inherent mathematical complexity and associated implementation challenges have resulted in a relative lack of attention to this concept in the literature.

In this paper, we propose a novel variationally consistent membrane wrinkling model based on the spectral decomposition of the strain tensor. We split the strain tensor into positive and negative components based on their eigenvalues. It allows us to decompose the strain energy density additively into positive and negative contributions. Considering only the positive part of the strain energy in the variational formulation, we obtain a membrane with zero compressive stiffness. By adding the negative strain energy multiplied with a very small factor, we further obtain a residual compressive stiffness, which improves stability and allows handling also states of slackening. Thereby, the strain-based wrinkling criterion can be self-adaptively satisfied without the reliance on if-else logical judgment in the implementation. It ensures the consistency of the wrinkling model and avoids many iterations triggered by a membrane state transition. The spectral decomposition also provides information on the wrinkling direction. We consistently derive the modified strain energy density with respect to the strain variable to determine the new stress and material tensor formulations. Finally, we incorporate our proposed wrinkling model into an isogeometric membrane element and verify its effectiveness with benchmark problems.

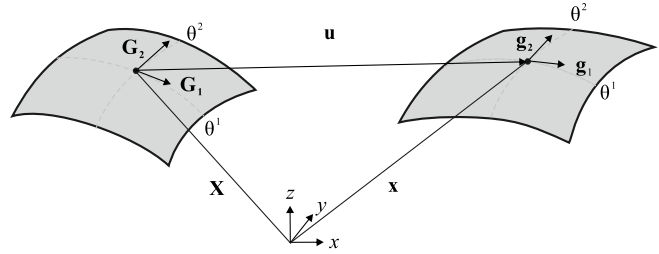

**Fig. 2.** Schematic description of a membrane: $\mathbf{X}$ and $\mathbf{x}$ are the position vectors on the midsurface in the reference configuration and the deformed configuration with displacement vector $\mathbf{u}$, respectively.

The paper is organized as follows. In Section 2, we provide an overview of the geometrical basics of the membrane description, the membrane kinematics, and its variational formulation. This section serves as the foundation for the subsequent sections. Section 3 presents the newly developed wrinkling model, which includes the spectral decomposition of the strain tensor and the consistent derivation of new stress and constitutive tensors. In addition, we briefly discuss the linearization of its variational formulation and isogeometric discretization in Section 4. In order to demonstrate the applicability of the presented methods, four tests are presented in Section 5, including benchmark examples from the literature with analytical, experimental, or numerical reference solutions. Finally, in Section 6, we summarize the key features of our proposed wrinkling model and suggest potential future directions for research in this area.

## 2. Membrane formulation

In this section, we present a concise overview of the fundamental aspects of the geometry, kinematics, and material law of a membrane. As shown in Fig. 2, a membrane can be generally represented as a curved surface with a certain thickness $t$, described by the position vector $\mathbf{x}(\theta^\alpha)$. The surface coordinates are denoted as $\theta^\alpha$. Greek indices take on values of $\{1, 2\}$, and the summation convention for repeated indices is employed. Then, the covariant base vectors $\mathbf{g}_\alpha$ on the midsurface in the current configuration are defined as:

$$\mathbf{g}_\alpha = \frac{\partial \mathbf{x}}{\partial \theta^\alpha} = \mathbf{x}_{,\alpha}. \tag{1}$$

Similarly, in the context of the reference configuration, the covariant basis vectors are denoted as $\mathbf{G}_\alpha$. Then, the contravariant basis vectors $\mathbf{G}^\alpha$ and $\mathbf{g}^\alpha$ can be determined by the rule:

$$\mathbf{G}^\alpha \cdot \mathbf{G}_\beta = \delta^\alpha_\beta, \quad \mathbf{g}^\alpha \cdot \mathbf{g}_\beta = \delta^\alpha_\beta, \tag{2}$$

where $\delta^\alpha_\beta$ is the Kronecker delta. If $\alpha = \beta$, $\delta^\alpha_\beta = 1$; otherwise, $\delta^\alpha_\beta = 0$. The metric coefficients of the midsurface in the reference configuration $G_{\alpha\beta}$ and current configuration $g_{\alpha\beta}$ are represented as:

$$G_{\alpha\beta} = \mathbf{G}_\alpha \cdot \mathbf{G}_\beta, \quad g_{\alpha\beta} = \mathbf{g}_\alpha \cdot \mathbf{g}_\beta. \tag{3}$$

Subsequently, the deformation gradient tensor $\mathbf{F}$ can be introduced and defined as:

$$\mathbf{F} = \mathbf{g}_\alpha \otimes \mathbf{G}^\alpha, \quad \mathbf{F}^T = \mathbf{G}_\alpha \otimes \mathbf{g}^\alpha. \tag{4}$$

In order to characterize the nonlinear relationship between deformations and strains, the Green–Lagrange strain tensor is employed and given by:

$$\mathbf{E} = \frac{1}{2} \left( \mathbf{F}^T \mathbf{F} - \mathbf{I} \right), \tag{5}$$

where $\mathbf{I}$ denotes the identity tensor. The in-plane coefficients of $\mathbf{E}$ are obtained by the metric coefficients in the deformed and undeformed configurations as:

$$E_{\alpha\beta} = \frac{1}{2} \left( g_{\alpha\beta} - G_{\alpha\beta} \right). \tag{6}$$

The second Piola–Kirchhoff (PK2) stress tensor $\mathbf{S}$ is introduced as the energetically conjugate quantity to the Green–Lagrange strain tensor $\mathbf{E}$. It is derived from the strain energy density $\psi$ with respect to the strain $\mathbf{E}$. Here, the strain energy density of the St. Venant-Kirchhoff constitutive model with plane stress is given by:

$$\psi\left(\mathbf{E}\right) = \frac{\lambda}{2}\left(\mathrm{tr}\left(\mathbf{E}\right)\right)^2 + \mu\mathrm{tr}\left(\mathbf{E}^2\right) - \frac{\lambda^2}{2\left(\lambda + 2\mu\right)}\left(\mathrm{tr}\left(\mathbf{E}\right)\right)^2, \tag{7}$$

where $\lambda$ is the first Lamé parameter and $\mu$ the second Lamé parameter, respectively. The stress tensor $\mathbf{S}$ can be determined via the tangent material tensor $\mathbb{C}$ and its coefficients $S^{\alpha\beta}$ are then expressed in terms of the material tensor coefficients $\mathbb{C}^{\alpha\beta\gamma\delta}$ as:

$$S^{\alpha\beta} = \mathbb{C}^{\alpha\beta\gamma\delta} E_{\gamma\delta}. \tag{8}$$

By approximating the differential volume $dV$ as the product of the thickness $t$ of membrane and its midsurface differential area $dA$, i.e., $dV \approx t\,dA$, the virtual work $\delta W(\mathbf{u}, \delta\mathbf{u})$ can be formulated as follows:

$$\delta W\left(\mathbf{u}, \delta\mathbf{u}\right) = \delta W^{\mathrm{int}} - \delta W^{\mathrm{ext}} = \int_A \mathbf{S} : \delta\mathbf{E}\; t\; dA - \int_A \mathbf{f} \cdot \delta\mathbf{u}\; dA, \tag{9}$$

which involves the internal virtual work $\delta W^{\mathrm{int}}$ and external virtual work $\delta W^{\mathrm{ext}}$. Within the Eq. (9), $\delta\mathbf{E}$ represents the virtual strain, $\mathbf{f}$ denotes the external force, and $\delta\mathbf{u}$ represents the virtual displacement.

## 3. New wrinkling model based on spectral decomposition

In this section, we present the new wrinkling model. The fundamental concept is to split the strain tensor $\mathbf{E}$ into positive and negative components via a spectral decomposition

$$\mathbf{E} = E_\alpha \mathbf{n}_\alpha \otimes \mathbf{n}_\alpha, \tag{10}$$

where $E_\alpha$ represents the eigenvalues (principal strains) corresponding to the eigenvectors $\mathbf{n}_\alpha$ (principal directions). The positive and negative strain tensors, denoted as $\mathbf{E}^+$ and $\mathbf{E}^-$, respectively, are then determined in terms of the positive and negative principal strains $\langle E_\alpha \rangle^\pm$

$$\mathbf{E}^\pm = \langle E_\alpha \rangle^\pm \mathbf{n}_\alpha \otimes \mathbf{n}_\alpha, \tag{11}$$

with $\langle x \rangle^\pm = (x \pm |x|)/2$, and the complete strain tensor is simply obtained by the sum of the positive and negative parts

$$\mathbf{E} = \mathbf{E}^+ + \mathbf{E}^-. \tag{12}$$

With the decomposition of the strain tensor $\mathbf{E}$, the strain energy density $\psi(\mathbf{E})$ can also be decomposed into positive and negative terms. This concept was presented first in [45] in the context of phase-field fracture and was extended in [46] to the case of plane stress problems as follows

$$\psi^\pm\left(\mathbf{E}\right) = \frac{\lambda}{2}\left(\langle\mathrm{tr}\left(\mathbf{E}\right)\rangle^\pm\right)^2 + \mu\mathrm{tr}\left(\left(\mathbf{E}^\pm\right)^2\right) - \frac{\lambda^2}{2\left(\lambda + 2\mu\right)}\left(\langle\mathrm{tr}\left(\mathbf{E}\right)\rangle^\pm\right)^2, \tag{13}$$

$$\psi\left(\mathbf{E}\right) = \psi^+\left(\mathbf{E}\right) + \psi^-\left(\mathbf{E}\right). \tag{14}$$

From Eq. (13) we can derive in a variationally consistent manner the split of the stress tensor into positive and negative terms

$$\mathbf{S}^\pm = \frac{\partial\psi^\pm\left(\mathbf{E}\right)}{\partial\mathbf{E}} = \left(\lambda - \frac{\lambda^2}{\lambda + 2\mu}\right)\langle\mathrm{tr}\left(\mathbf{E}\right)\rangle^\pm\mathbf{I} + 2\mu\mathbf{E}^\pm, \tag{15}$$

$$\mathbf{S} = \mathbf{S}^+ + \mathbf{S}^-. \tag{16}$$

In the same way, we split the constitutive tensor in positive and negative terms

$$\mathbb{C}^\pm = \frac{\partial\mathbf{S}^\pm\left(\mathbf{E}\right)}{\partial\mathbf{E}} = \left(\lambda - \frac{\lambda^2}{\lambda + 2\mu}\right)H\left(\pm\mathrm{tr}\left(\mathbf{E}\right)\right)\mathbb{J} + 2\mu\left(H\left(\pm E_\alpha\right)\mathbb{Q}_\alpha + \sum_{\alpha\neq\beta}^2 \frac{\langle E_\alpha \rangle^\pm}{2\left(E_\alpha - E_\beta\right)}\left(\mathbb{G}_{\alpha\beta} + \mathbb{G}_{\beta\alpha}\right)\right), \tag{17}$$

$$\mathbb{C} = \mathbb{C}^+ + \mathbb{C}^-. \tag{18}$$

with the fourth-order tensor $(\mathbb{Q}_\alpha)_{\gamma\delta\epsilon\zeta} := (\mathbf{M}_\alpha)_{\gamma\delta}(\mathbf{M}_\alpha)_{\gamma\delta}$. The detailed derivation of Eq. (17) is presented in  Appendix A.

In Eqs. (11)–(17), positive terms may be associated to tension and negative terms to compression. A membrane formulation according to the assumption of tension field theory, in which only tensile stresses are possible, can then easily be obtained by only considering the positive parts in Eqs. (14)–(18)

$$\psi = \psi^+, \tag{19}$$

$$\mathbf{S} = \mathbf{S}^+, \tag{20}$$

$$\mathbb{C} = \mathbb{C}^+, \tag{21}$$

and a finite element (or isogeometric) formulation can be derived in a standard way as shown in Section 4. Notably, no explicit wrinkling criterion needs to be evaluated at the material point, and no explicit distinction between taut, wrinkled, or slack, with different formulations for the different cases, is necessary during analysis. Eqs. (19)–(21) include all three cases. Nevertheless, the model also provides the relevant information about wrinkling, i.e., the membrane state (taut, wrinkled, or slack) and the wrinkling direction. These are obtained automatically when doing the spectral decomposition of the strain tensor in Eq. (10) as follows

$$E_1 > 0, \quad \begin{matrix} E_2 > 0 : & \text{taut}, \\ E_2 \leq 0 : & \text{wrinkled}, \\ E_1 \leq 0 : & \text{slackened}, \end{matrix} \tag{22}$$

with $E_1 \geq E_2$, and, in the case of wrinkling, the wrinkling direction is obtained as $\mathbf{n}_1$, the eigenvector corresponding to $E_1$. The formulation based on Eqs. (19)–(21) describes a membrane with zero compressive stiffness, which is in accordance with tension field theory, allowing only for taut and wrinkled states, but not for slack. In a corresponding finite element formulation, even slight compression may lead to self-penetration of elements and slack states, even if they occur only locally, may lead to singular stiffness matrices. An efficient way to prevent this is to assign some residual compressive stiffness. In the presented approach, this can be done simply by modifying Eqs. (19)–(21) as follows

$$\psi = \psi^+ + \eta \psi^-, \tag{23}$$

$$\mathbf{S} = \mathbf{S}^+ + \eta \mathbf{S}^-, \tag{24}$$

$$\mathbb{C} = \mathbb{C}^+ + \eta \mathbb{C}^-, \tag{25}$$

where $\eta \ll 1$ is a factor governing the compressive stiffness. Theoretically, this should be a function of the buckling strength of the membrane. In this work, however, we use it as a prescribed stabilization parameter, which helps to prevent self-penetration and singular stiffness matrices, and, moreover, can improve significantly the iterative convergence. The idea of assigning a small compressive stiffness via a prescribed parameter can also be found in other works in the literature, e.g., in [2,18], where the parameter is considered as a penalty parameter, in [38], where compressive stresses are termed as "allowable compressive stress", or in [1], where compressive stiffness is termed as "residual compressive stiffness". In this paper, we adopt the terminology of "residual" compressive stiffness and stresses, which are governed by the "degradation factor" $\eta$. In Section 5.3, we show a numerical example where different values for $\eta$ and its impact on the results are studied. For many problems, where slacking or self-penetration is not relevant, $\eta = 0$ can be used.

## 4. Discretization

In this section, we present the linearization of Eq. (9) and discretization procedures for solving this equation. The presented discretization is suitable for any basis function, such as the Lagrange polynomial used in standard finite element analysis or Non-Uniform Rational B-Splines (NURBS) typically applied in isogeometric analysis (IGA), to name a few. IGA is a widely adopted technique in structural analysis due to its numerous advantages. To ensure conciseness, we refrain from reiterating these details here; more details can be found in the literature [47,48].

The discretized displacement field is represented as follows:

$$\mathbf{u} = \sum_{a}^{n_{sh}} N^a \mathbf{u}^a, \tag{26}$$

where $N^a$ denotes the shape functions, $n_{sh}$ represents the total number of shape functions, and $\mathbf{u}^a$ refers to the nodal displacement vectors with components $u_i^a$ ($i = 1, 2, 3$) representing the global $x$-, $y$-, $z$-components. We establish an expression for the global degree of freedom number $r$ of a nodal displacement, i.e., $r = 3(a-1) + i$, such that $u_r = u_i^a$. In order to determine the variation of the displacement field with respect to $u_r$, the partial derivative of $\partial / \partial u_r$ is used:

$$\frac{\partial \mathbf{u}}{\partial u_r} = N^a \mathbf{e}_i, \tag{27}$$

where $\mathbf{e}_i$ denote the global Cartesian base vectors. For more details, we refer to Kiendl et al. [49]. Subsequently, by deriving the variations of the internal virtual work $\delta W^{\text{int}}$ and external virtual work $\delta W^{\text{ext}}$ in Eq. (9) with respect to $u_r$, the residual force vector $\mathbf{R}$ is obtained and defined as:

$$R_r = F_r^{\text{int}} - F_r^{\text{ext}} = \int_A \mathbf{S} : \frac{\partial \mathbf{E}}{\partial u_r} \, t \, \mathrm{d}A - \int_A \mathbf{f} \cdot \frac{\partial \mathbf{u}}{\partial u_r} \, \mathrm{d}A, \tag{28}$$

where $\mathbf{F}^{\text{int}}$ and $\mathbf{F}^{\text{ext}}$ represent the vectors of the internal and external nodal loads, respectively. The linearization of Eq. (28) yields the tangential stiffness matrix $\mathbf{K}$, which comprises the internal stiffness matrix $\mathbf{K}^{\text{int}}$ and external stiffness matrix $\mathbf{K}^{\text{ext}}$, formulated as follows:

$$K_{rs} = K_{rs}^{\text{int}} - K_{rs}^{\text{ext}} = \int_A \frac{\partial \mathbf{S}}{\partial u_s} : \frac{\partial \mathbf{E}}{\partial u_r} + \mathbf{S} : \frac{\partial^2 \mathbf{E}}{\partial u_r \partial u_s} \, t \, \mathrm{d}A - \int_A \frac{\partial \mathbf{f}}{\partial u_s} \cdot \frac{\partial \mathbf{u}}{\partial u_r} \, \mathrm{d}A. \tag{29}$$

In order to solve the linearized equation system, the Newton–Raphson method is employed and given by:

$$\frac{\partial W}{\partial u_r} + \frac{\partial^2 W}{\partial u_r \partial u_s} \Delta u_s = 0, \tag{30}$$

where $\Delta u_s$ denotes the components of the incremental displacements. Through the solution of the equation system, as mentioned earlier, we derive the incremental displacement vector $\Delta \mathbf{u}$, which ensures its accuracy subject to the precise computation of the residual vector $\mathbf{R}$, given by:

$$\mathbf{K}\Delta \mathbf{u} = -\mathbf{R}. \tag{31}$$

The Newton–Raphson method converges towards the desired solution by iteratively updating the displacement vector until the residual is minimized, leading to the accurate determination of the incremental displacement $\Delta \mathbf{u}$.

## 5. Numerical examples

In this section, we assess the capability of the proposed formulation in four numerical examples with analytical or reference solutions from the literature. All tests are conducted using nonlinear static isogeometric analysis, performed with NURBS. In order to show that our proposed formulation is applicable to standard finite element analysis, we also set the polynomial degree of the NURBS basis function as $p = 1$ for some cases because it is equivalent to first-order Lagrange polynomial when $p = 1$. The following notations are used for various stresses presented in the results: $\sigma$ indicates the Cauchy stress tensor, which is obtained by $\sigma = (\det \mathbf{F})^{-1} \cdot \mathbf{F} \cdot \mathbf{S} \cdot \mathbf{F}^T$. Its first and second principal stresses are denoted as $\sigma_1$ and $\sigma_2$, respectively. The residual norm applied in the analysis is defined with respect to the external forces magnitude as $\|\mathbf{R}\|/\|\mathbf{F}^{\text{ext}}\|$. Unless stated otherwise, a degradation factor $\eta = 0$ is used.

### 5.1. In-plane pure bending of a pre-tensioned rectangular membrane

The first numerical example considers a pre-stretched rectangular membrane under in-plane pure bending. Since an analytical solution was presented in [5], this example is considered a benchmark for studying partly wrinkled membranes. This problem has also been extensively investigated in the literature by various researchers as [8,27,38,50–53]. As depicted in Fig. 3, the membrane with the height $H$ and the thickness $t$ is subjected to a uniform stress $\sigma_0$ in the $y$-direction, as well as a pair of the axial loads $P = \sigma_0 t H$ and the bending moments $M$ at the lateral sides. By increasing the bending moments $M$, a band of vertical wrinkles of height $h$ (highlighted in gray) appears along the bottom edge, and the band height $h$ of the wrinkled zone is defined in [5] as:

$$\frac{h}{H} = \begin{cases} 0 & M/PH < 1/6 \\ 3M/PH - 1/2 & 1/6 \le M/PH < 1/2. \end{cases} \tag{32}$$

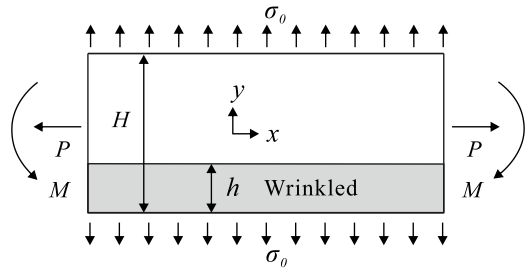

**Fig. 3.** In-plane pure bending of a pre-tensioned rectangular membrane.

In the wrinkled region, the normal stress $\sigma_x$ in the membrane is eliminated, whereas, in other regions, it is linearly distributed along the height $H$. The normal stress distribution can be calculated using the following expression:

$$\frac{\sigma_x}{\sigma_0} = \begin{cases} \frac{2(y/H - h/H)}{(1 - h/H)^2} & h/H < y/H \le 1 \\ 0 & 0 \le y/H \le h/H. \end{cases} \tag{33}$$

It must be noted that this analytical solution relies on a stress-based wrinkling criterion. Accordingly, the proposed model is expected to converge to this solution only for Poisson's ratio $\nu = 0$. To verify this, we perform simulations with $\nu = 0$ and $\nu = 0.3$ in the following.

The simulation setup for this problem is presented in Fig. 4, where we use 55 ($11 \times 5$) bi-quadratic isogeometric membrane elements to solve the problem, following the recommendation in [29]. Only the right half of the membrane is modeled due to its symmetry, and the displacements in the $x$-direction on the left are constrained. The middle point of the left edge is also prevented from deforming in the $y$-direction. As illustrated, the axial load $P$ and bending moment $M$ are replaced by equivalent stresses $\sigma_p = P/tH$ and $\sigma_M = 6M/tH^2(2y/H - 1)$, respectively. In order to withstand compressive stresses and maintain a uniform rotation, the five elements on the right are assumed to remain taut. Thus, these elements are modeled as the standard membrane elements without embedding the wrinkling model.

Fig. 5 shows the normal stress $\sigma_x$ distribution with respect to the applied stress $\sigma_0$ along the height $H$. Here, we use $y/H$-ratio, i.e., the stress-measured position $y$ related to the height $H$, as $x$-axis to unify the results. Different curves exhibit the results for different values of the bending moment-force ratio (2M/PH). For the case of Poisson's ratio $\nu = 0$, the results of the proposed model coincide with the analytical solution, as shown in Fig. 5(a), while a slight difference is observed for the case of $\nu = 0.3$, as shown in Fig. 5(b). This deviation is caused by the fact, as mentioned above, that the proposed model relies on a strain-based wrinkling criterion, the analytical solution on a stress-based criterion. Thus, the second principal stress $\sigma_2$ evaluated by our model is not always equal to zero. In contrast, the analytically computed stresses are enforced to satisfy the uniaxial tension condition. Although the Poisson effect could influence the stress response predicted by the proposed model compared to the analytical solutions, the difference observed in our model can be neglected in practical applications, as demonstrated by the following numerical examples.

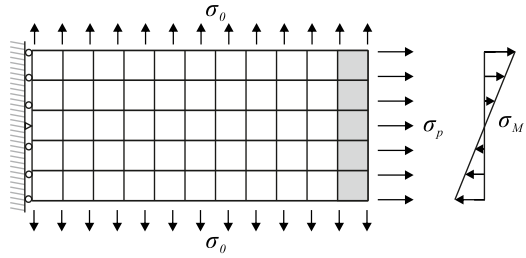

**Fig. 4.** Numerical simulation setting for the right half of the rectangular membrane under in-plane bending.

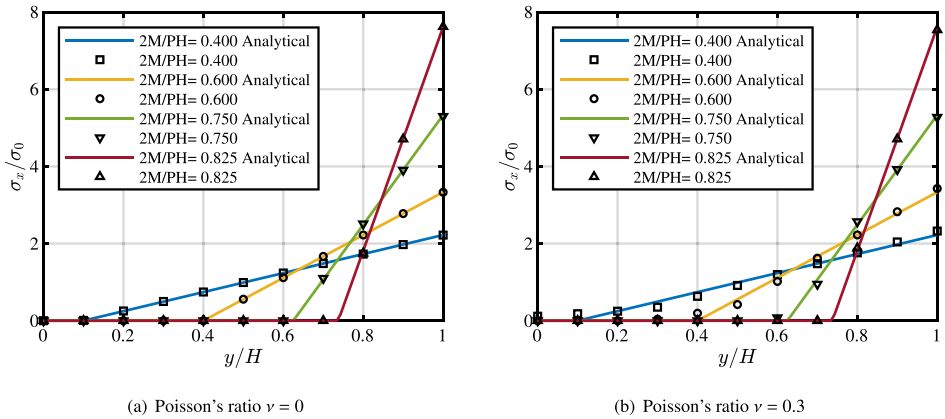

(a) Poisson's ratio $\nu = 0$

(b) Poisson's ratio $\nu = 0.3$

**Fig. 5.** Normal stress $\sigma_x$ distribution along the height $H$ with the newly proposed wrinkling model.

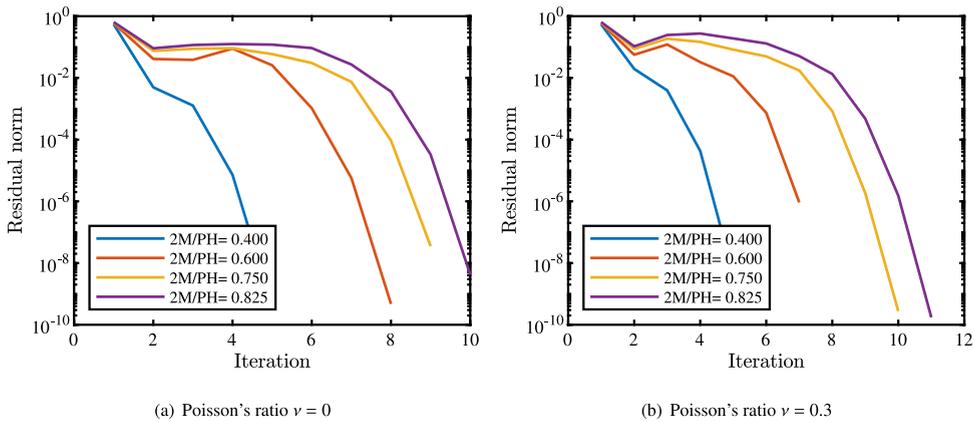

(a) Poisson's ratio $\nu = 0$

(b) Poisson's ratio $\nu = 0.3$

**Fig. 6.** Convergence behavior of rectangular membrane under in-plane bending during the second load step.

In addition, in order to show the convergence behavior of the proposed model, we plot in Fig. 6 the number of iterations during the second load step, where we can find the most iterations of all load steps. From Figs. 6(a) and 6(b), it can be seen that the new model can achieve a good convergence and increasing bending moment requires more iterations to reach the preset tolerance. To prove the robustness of the proposed wrinkling model during mesh refinement, we conduct a mesh convergence study by plotting the normal stresses $\sigma_x$ along with the height $H$ under various mesh sizes with a Poisson's ratio $\nu = 0$ in Fig. 7. The results of the mesh convergence studies confirm that the accuracy and robustness of the proposed model are not affected by changes in the element size. In Appendix B, the results for $p = 1$ are attached in Fig. B.25 to evidence that the newly proposed model can be equally applied to standard finite element analysis.

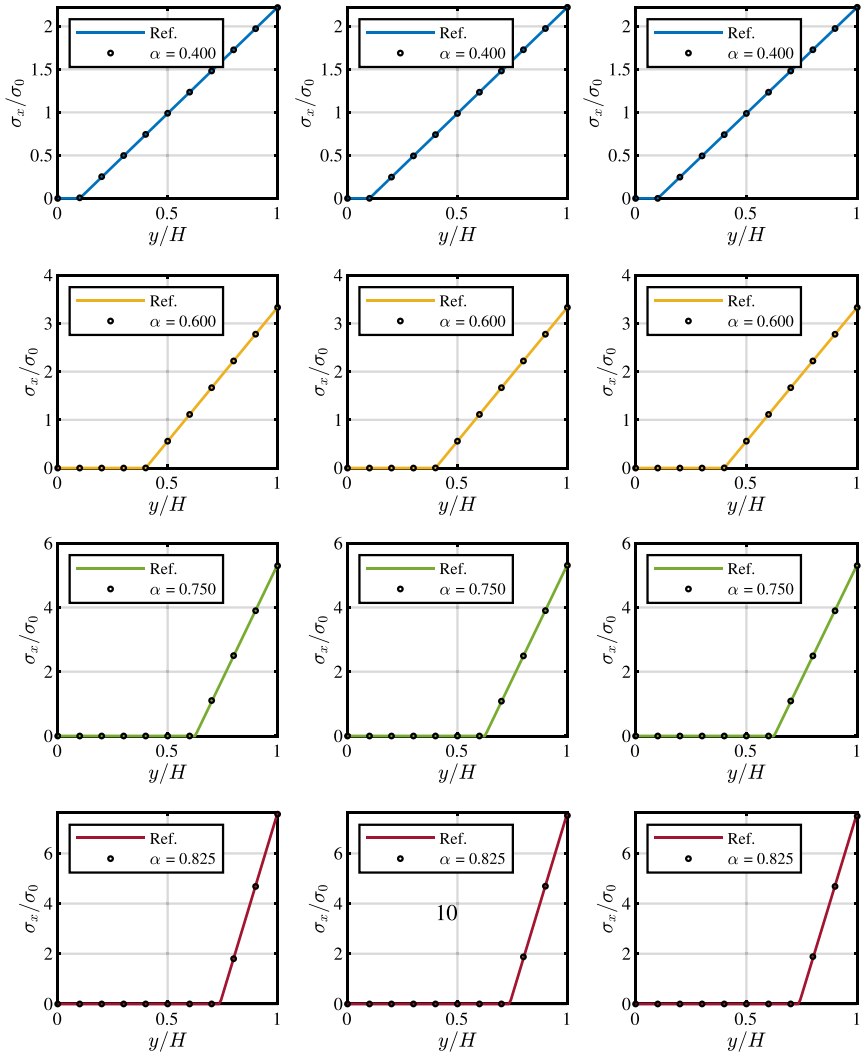

**Fig. 7.** Normal stress $\sigma_x$ distribution along the height $H$ computed by the proposed model with Poisson's ratio $\nu = 0$, polynomial degree $p = 2$, and $\alpha = 2M/PH$. In the left, middle, and right columns, mesh refinement is carried out as $11 \times 5$, $23 \times 11$, and $43 \times 21$, respectively.

### 5.2. Flat square membrane under corner loads

As a second example, we present a verification of the proposed wrinkling model using the case of a square membrane subjected to corner loads of $T_1$ and $T_2$ with nonzero Poisson's ratio, which is discussed in [4,16,38,54]. The membrane possesses a length of $0.5$ m and a thickness of $25\,\mu\text{m}$. It undergoes the diagonal pairs of equal and opposite forces $T_1$ and $T_2$ at the four corners, as shown in Fig. 8. As suggested in [38], the central point of the membrane is fixed to prevent rigid body motions. In addition, the movement in $x$-direction of the middle point at the top edge is also constrained to avoid rotation. We perform numerical analyses for the various values of the $T_1/T_2$-ratio, where $T_1$ is increased from $5\,\text{N}$ to $20\,\text{N}$, while $T_2$ is kept constant at $5\,\text{N}$ for all cases. The material parameters used in this study are Young's modulus $E = 3500\,\text{MPa}$ and Poisson's ratio $\nu = 0.31$. The square membrane is discretized using $40 \times 40$ cubic elements.

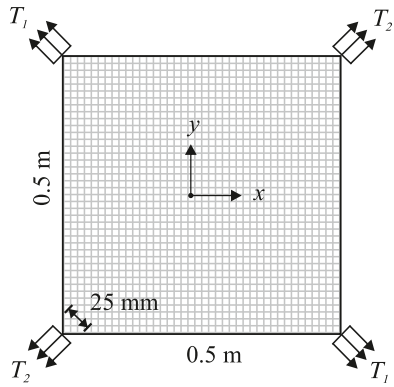

**Fig. 8.** Flat square membrane under corner loads $T_1$ and $T_2$.

We plot the first principal stress distributions $\sigma_1$ under various load ratios $T_1/T_2$ in Fig. 9(a)–(d). The new model appears to attain the comparable first principal stress distribution patterns under various load ratios $T_1/T_2$ compared to the thin shell solutions [16] and those predicted by other wrinkling models [4,38]. Additionally, we plot the wrinkling intensity under the different values of the $T_1/T_2$ ratio in Fig. 10. The wrinkling intensity is obtained as the absolute values of the negative second eigenvalues of the strain tensor multiplied with the norm of the first eigenvector if the wrinkling state is wrinkled. Otherwise, in the case of taut or slackened states, the value would be zero. In Fig. 11, the wrinkling trajectories are presented. The regions marked by lines indicate the wrinkled zones, the lines show the direction of the first eigenvector of the strains, and the line length is scaled according to the corresponding second eigenvalues. The areas without lines are the taut or slackened regions. As shown in Fig. 11, it can be recognized that the wrinkles are mainly concentrated in the diagonal by increasing the load ratio, similar to the fact observed in the experiments [54] and the numerical results [4,16,38]. Fig. 12 presents a graphical representation illustrating the convergence performance of the square membrane subjected to various corner load ratios. The convergence performance is measured by tracking the number of iterations required in the first and last load steps. As depicted in Fig. 12(a), in the case of ratio $T_1/T_2 = 1$, it requires fewer iterations to reach the tolerance compared to other cases because the diagonal wrinkles that can deteriorate the convergence behavior have not formed yet in the case of $T_1/T_2 = 1$. The convergence behavior from the second load step follows a similar pattern to the last load step. Hence, we show the iterations in the last load step as illustrated in Fig. 12(b).

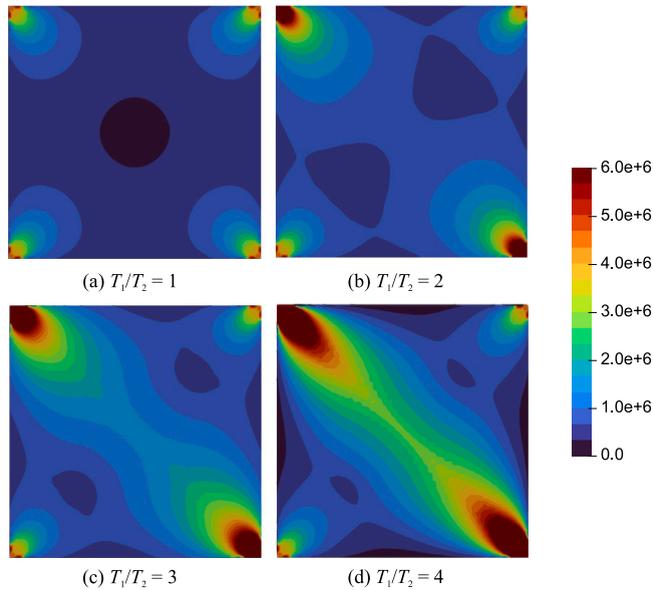

(a) $T_1/T_2 = 1$

(b) $T_1/T_2 = 2$

(c) $T_1/T_2 = 3$

(d) $T_1/T_2 = 4$

**Fig. 9.** Contours of the first principal stresses $\sigma_1$ (Pa) under the different values of the $T_1/T_2$-ratio.

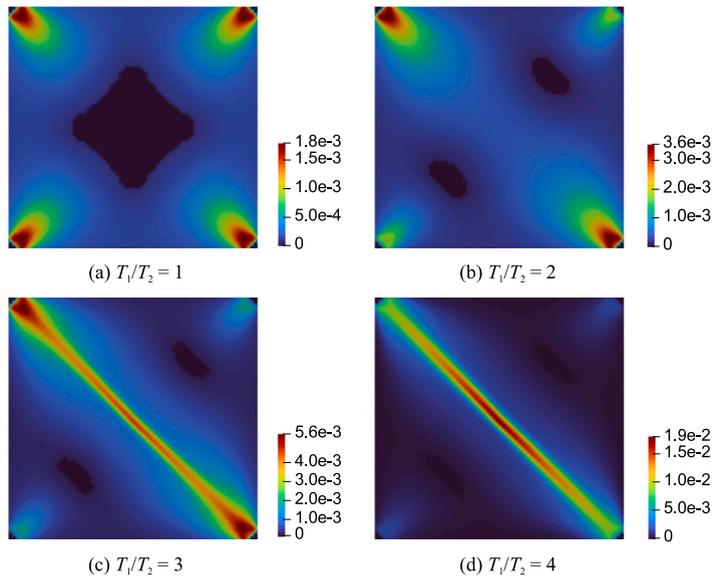

**Fig. 10.** Wrinkling intensity under the different values of the $T_1/T_2$-ratio.

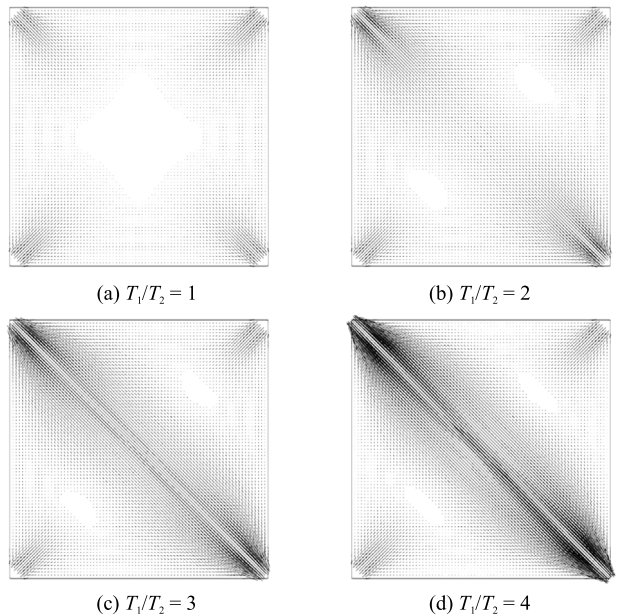

**Fig. 11.** Wrinkling trajectories under the different values of the $T_1/T_2$-ratio.

### 5.3. Inflation of a square isotropic airbag

In the third example, we model the inflation of a square airbag, which is a commonly used benchmark for validating finite element-based wrinkling algorithms; related works are referred to Diaby et al. [55], Jarasjarungkiat et al. [38], Jarasjarungkiat et al. [18], Contri and Schrefler [6], Kang and Im [34], Lee and Youn [56], Le Meitour et al. [57] and Gil and Bonet [58]. The initially flat square airbag has a diagonal length of $AC = 120$ cm and a thickness of $t = 0.06$ cm, as Fig. 13 shows. It is loaded by a displacement-dependent pressure $P$ perpendicular to the surface, increasing to $5000$ Pa gradually. We did not consider the external stiffness matrix produced by pressure in our element formulation for simplification. Due to the symmetry, only a quarter of the airbag is simulated. Thus, the symmetric boundary conditions are applied to the inner edges, while the movements of outer edges in the $z$−direction are constrained. The material of the airbag is assumed to be linear isotropic, with an elastic modulus $E = 588$ MPa and Poisson's ratio $\nu = 0.4$.

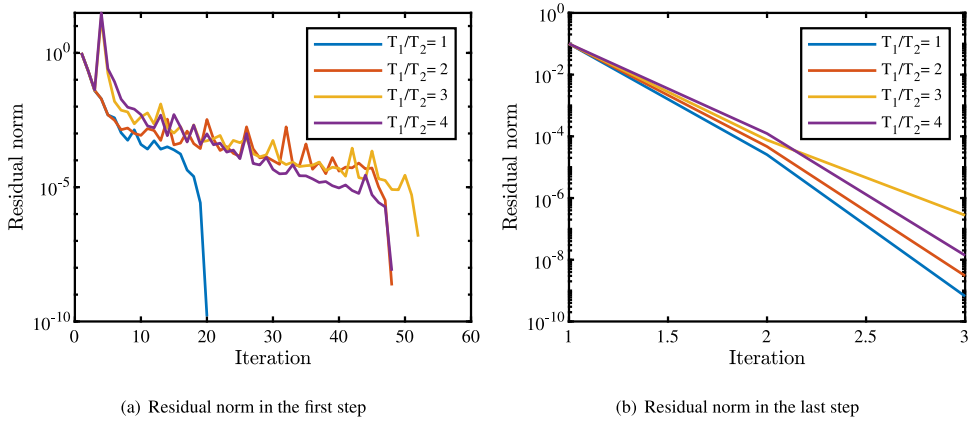

(a) Residual norm in the first step

(b) Residual norm in the last step

**Fig. 12.** Convergence behavior of the square membrane under corner loads with various $T_1/T_2$-ratios.

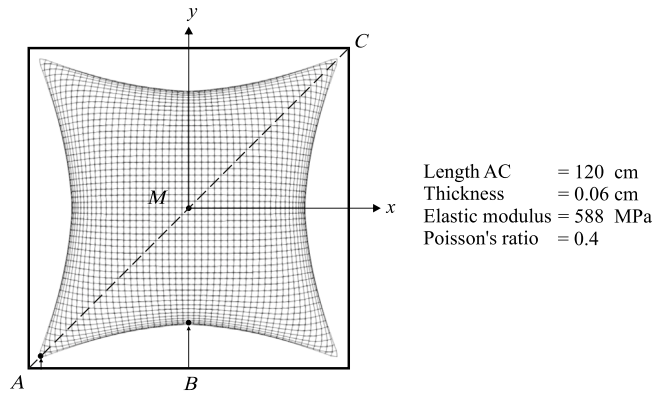

| | |
|---|---|
| Length AC | = 120 cm |
| Thickness | = 0.06 cm |
| Elastic modulus | = 588 MPa |
| Poisson's ratio | = 0.4 |

**Fig. 13.** Airbag simulation setup.

As indicated in [55], the inflation of a square airbag poses a specific difficulty caused by the singularity of the stiffness matrix at the initial stage due to the normal pressure applied to the membrane surface. One possibility is to use dynamic relaxation (DR) [32] to overcome this numerical obstacle. DR approximates a static equilibrium solution through a pseudo-dynamic transient analysis. However, a more straightforward alternative is to stretch the membrane along the $x$- and $y$-directions with dead forces on the two outer edges, as suggested in [55,56]. These forces are gradually reduced and removed when the pressure $P$ reaches a fixed value. This approach enables us to avoid the initial singularity and obtain an accurate solution for the airbag inflation simulation.

In Table 1, we report the displacements and the first principal stress obtained with the proposed wrinkling model at various points, as the vertical displacement $u_z^M$ at point $M$ and the displacements of the $y$-components at points $A$ and $B$, denoted as $u_y^A$ and $u_y^B$, respectively. Moreover, the first principal stress $\sigma_1^M$ at point $M$ was also measured.

**Table 1**
Airbag results under various mesh size with $\eta = 0$.

| Number of elements | $2 \times 2$ | $4 \times 4$ | $8 \times 8$ | $16 \times 16$ | $32 \times 32$ |
|---|---|---|---|---|---|
| $u_z^M$ (m) | 0.2234 | 0.2175 | 0.2172 | 0.2170 | 0.2174 |
| $u_y^A$ (m) | 0.0396 | 0.0378 | 0.0368 | 0.0358 | 0.0352 |
| $u_y^B$ (m) | 0.1077 | 0.1168 | 0.1228 | 0.1289 | 0.1346 |
| $\sigma_1^M$ (MPa) | 13.91 | 3.95 | 4.27 | 3.90 | 3.80 |

Particularly noteworthy for discussion and analysis is the case where the mesh is $32 \times 32$. In this scenario, we observed a sudden change in the displacement $u_z^M$ at the center of the airbag. Similarly, a sudden contraction was observed at point $B$. This is attributed to the occurrence of self-penetration, as depicted in Fig. 14. This phenomenon has also been documented in [32]. Other

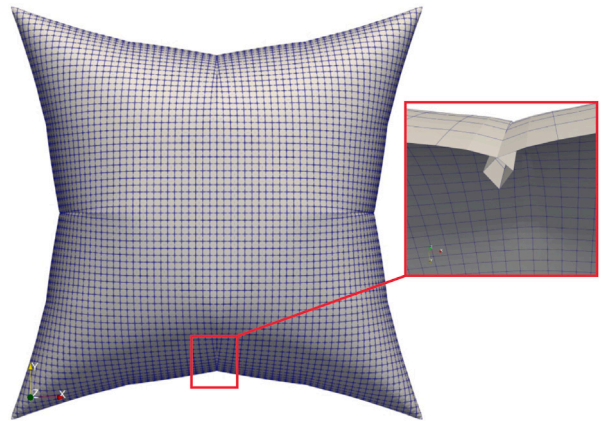

**Fig. 14.** Self-penetration at the point B of airbag with $32 \times 32$ mesh ($\eta = 0$).

published results did not perform a mesh convergence study and therefore did not refine the mesh to this extent. Consequently, in this example, penetration has been rarely discussed. Our interpretation of this is based on the studies of Iwasa et al. [19], Miyazaki [1], Jarasjarungkiat et al. [18,38], who pointed out that membranes under wrinkling conditions still should exhibit a small compressive stiffness, which is considered via the parameter $\eta$ in Eqs. (23)–(25) in this work. In the following, we perform a parametric study on $\eta$ and its influence on the results, where $\eta = 10^{-8}, 10^{-6}, 10^{-5}, 10^{-4}$ are considered. In Fig. 15, we plot $u_z^M$ and $u_y^B$ under different values of $\eta$ and increasingly refined meshes. It can be seen that for $\eta$ ranging from $10^{-8}$ to $10^{-5}$, no convergence for $u_z^M$ and $u_y^B$ is obtained within the meshes considered, while good convergence behavior is obtained with $\eta = 10^{-4}$. Fig. 16 shows exemplarily the detailed deformation at point *B* for the mesh with $32 \times 32$ elements and different values for $\eta$. It can be seen that for $\eta$ ranging from $10^{-8}$ to $10^{-5}$ (Fig. 16(a)–(c)), there is self-penetration, which disappears for $\eta = 10^{-4}$, see Fig. 16(d). This example shows that a

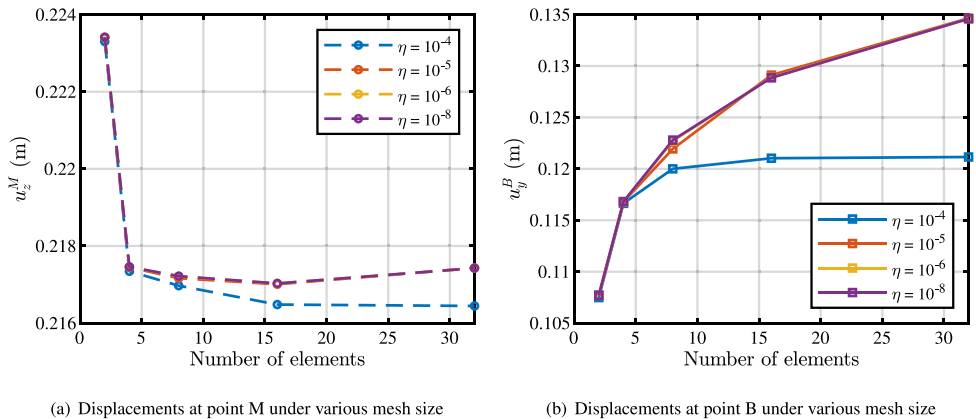

(a) Displacements at point M under various mesh size

(b) Displacements at point B under various mesh size

**Fig. 15.** Displacements at points M and B under various $\eta$.

small amount of compressive stiffness, governed via the parameter $\eta$ can avoid unphysical self-penetration and improve the overall convergence behavior. Regarding the choice of the value for $\eta$, the study shows that it may not have the desired effect if $\eta$ is chosen too low (see Fig. 16). At the same time, $\eta$ must not be chosen too high, as this would lead to a significant compressive stiffness in the membrane, affecting its global mechanical behavior. From Fig. 15 it becomes also clear that the choice of appropriate values for $\eta$ can depend on the mesh size. For the very coarse meshes ($2 \times 2$ and $4 \times 4$ elements), the choice of $\eta$ has no influence on the results, while for $8 \times 8$ elements and finer, it does. In conclusion, this study shows that the optimal value for $\eta$ can depend on the specific problem and the mesh size, and has to be determined empirically. At the same time, it has to be highlighted that even with $\eta = 0$ (or values of $\eta$ which are not optimal), the global behavior of the membrane is represented very well and also the local values for displacements and stresses are still in the range of reference results from the literature.

A comprehensive comparison between the results from the existing literature and the present work is presented in Table 2. By comparing our computed results with the outcomes reported in previous studies, a substantial agreement is observed, particularly concerning the prediction of the vertical displacements $u_z^M$ and the first principal stresses $\sigma_1^M$. It indicates that the slight difference due to the Poisson effect shown in the first benchmark can be neglected in complicated examples. In order to demonstrate the mesh independence of the proposed wrinkling formulation, we focus on analyzing the vertical displacements at point $M$ and the first principal stresses at the same location while progressively refining the mesh using linear ($p = 1$) and quadratic ($p = 2$) membrane elements. The results are presented in Fig. 17(a) for the vertical displacements and in Fig. 17(b) for the first principal stresses.

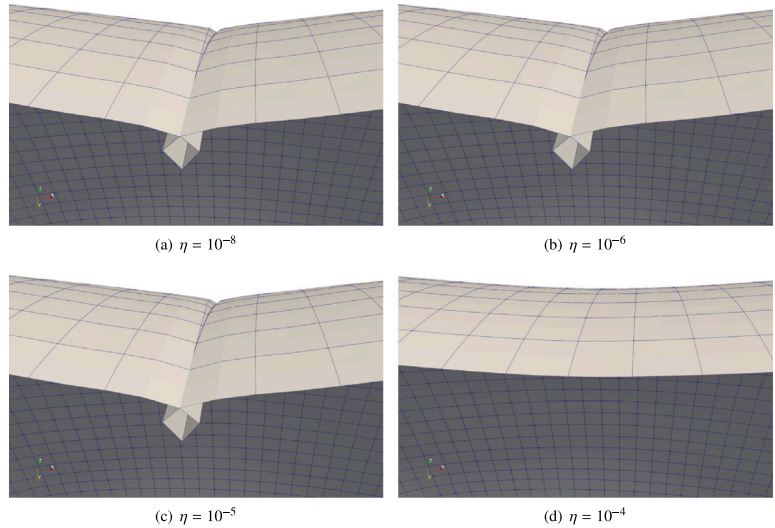

(a) $\eta = 10^{-8}$          (b) $\eta = 10^{-6}$

(c) $\eta = 10^{-5}$          (d) $\eta = 10^{-4}$

**Fig. 16.** Detail view at point $B$ for the $32 \times 32$ mesh with different values of $\eta$.

**Table 2**
Comparison of airbag results from literature and present work ($\eta = 10^{-4}$).

|  | Contri and Schrefler [6] | Kang and Im [34] | Diaby et al. [55] | Jarasjarungkiat et al. [38] | Present work $p = 1$ | Present work $p = 2$ |
|---|---|---|---|---|---|---|
| $u_z^M$ (m) | 0.217 | 0.214 | 0.2245 | 0.2175 | 0.2165 | 0.2164 |
| $u_x^A$ (m) | 0.045 | 0.041 | 0.0307 | 0.0349 | 0.0362 | 0.0351 |
| $u_y^B$ (m) | 0.110 | 0.119 | 0.1158 | 0.1203 | 0.1210 | 0.1212 |
| $\sigma_1^M$ (MPa) | 3.5 | – | – | 3.9 | 3.9 | 3.9 |

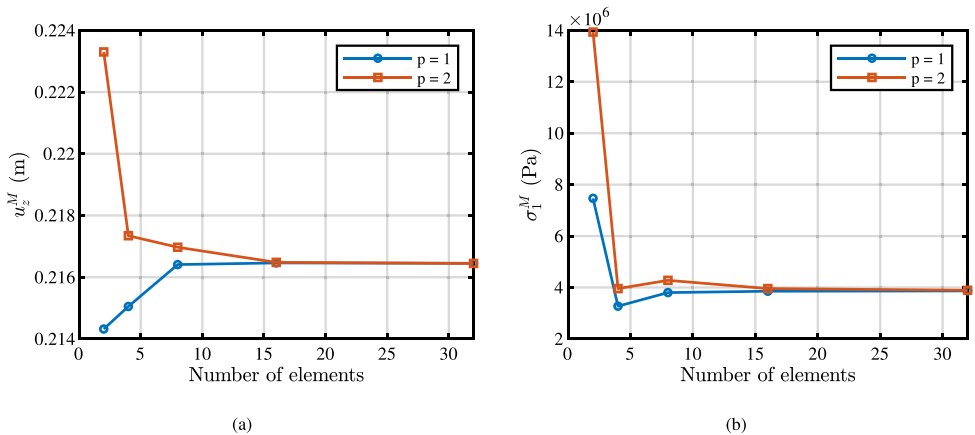

(a)                                    (b)

**Fig. 17.** Convergence study of the vertical displacement (a) and the first principal stress (b) at point $M$ during mesh refinement.

Besides, we plot the distributions of the first and second principal stresses $\sigma_{1/2}$ within the inflated airbag shown in Fig. 18. The maximum values of the first and second principal stresses are approximately $\sigma_1 = 49$ MPa and $\sigma_2 = 17$ MPa, both occurring at the corners. However, due to the excessive magnitude of these values, the stress distribution diagram lacks clarity. Therefore, in the stress contour plot, we present a more reasonable and scaled representation. Since the degradation factor is not set to zero, a small compressive stiffness remains within the model. Consequently, this leads to the generation of negative stresses that are small compared to the positive stresses and can be considered negligible. Besides, These principal stress distributions align with the observed wrinkle trajectories illustrated in Fig. 19. Precisely, the regions in proximity to the midpoint of the edges of the airbag exhibit pronounced wrinkling phenomena, while the central regions remain taut and devoid of wrinkles. Fig. 20 showcases the convergence behavior of the airbag during the final load step under different mesh sizes. The plot reveals a noteworthy trend wherein the mesh refinement leads to increased iterations required to attain the predetermined residual criterion. Furthermore, when the polynomial order is elevated, more iterations are necessary.

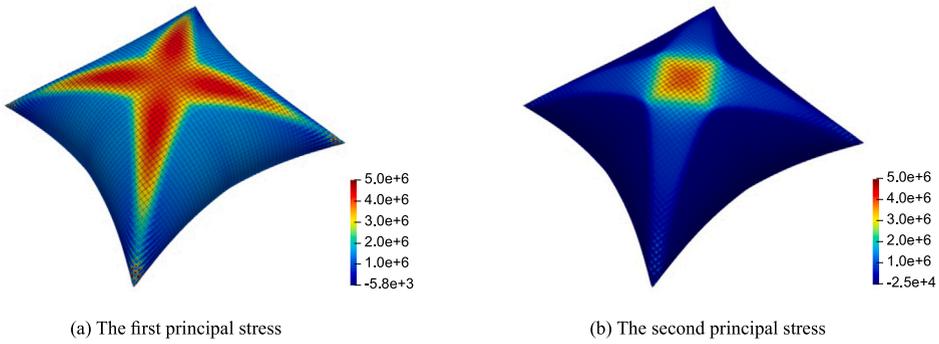

(a) The first principal stress                   (b) The second principal stress

**Fig. 18.** Contours of the first and second principal stresses $\sigma_{1/2}$ (Pa) distributed in the inflated airbag.

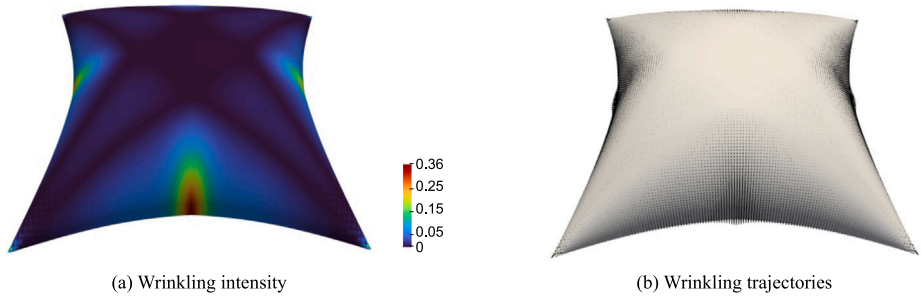

(a) Wrinkling intensity                       (b) Wrinkling trajectories

**Fig. 19.** Illustration of the wrinkled zone in a fully inflated square airbag.

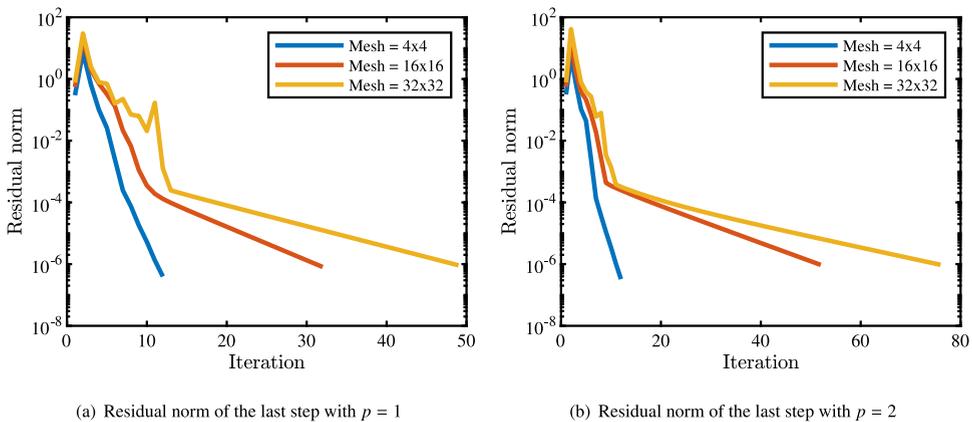

(a) Residual norm of the last step with $p = 1$           (b) Residual norm of the last step with $p = 2$

**Fig. 20.** Convergence behavior of airbag during the last step with different mesh.

### 5.4. Hanging blanket under self-weight

While all the benchmark examples above represent prestretched membranes, in this final example, we consider a "loose" membrane loaded by self-weight only. It represents a hanging blanket supported at its corners, as shown in Fig. 21. All four corners are supported rigidly in $z$-direction, while elastic supports are applied in the $x$-$y$-plane as shown in Fig. 21 such that the blanket can undergo large deformations due to self-weight. The length of the blanket is $L = 1$ m and its thickness $t = 1.177$ mm. The material parameters applied in this example are Youngs's modulus $E = 30\,000$ Pa, surface density $\rho = 0.144$ kg m$^{-2}$ for self-weight, which are taken from [59]. Poisson's ratio varies between $\nu = 0$ and $\nu = 0.3$, and the elastic supports are applied via the penalty formulation presented in [60], with a penalty parameter $\alpha = 10^{-2}$, corresponding to a spring stiffness $k_{spring} = 22.95$ kN m$^{-1}$. A mesh of $25 \times 25$ bi-quadratic elements is used to model the wrinkling behavior of the square hanging blanket.

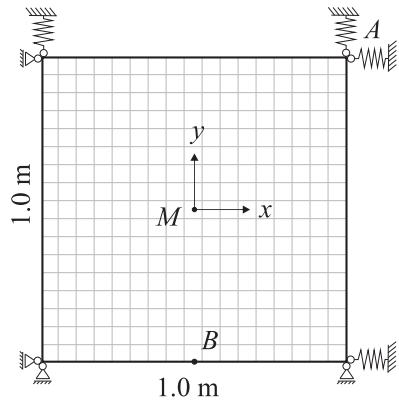

**Fig. 21.** Hanging blanket under self-weight.

Fig. 22 depicts the deformation and the first and second principal stresses. In Fig. 23, the wrinkling intensity and trajectories are plotted. As can be seen, significant wrinkling phenomena are prominently observed in the regions near the four corners and four edges of the blanket. In contrast, the central regions remain taut without any wrinkles. For validation, we also perform simulations with the wrinkling model proposed in [29] using the same setup and mesh, and its results are remarked as a reference.

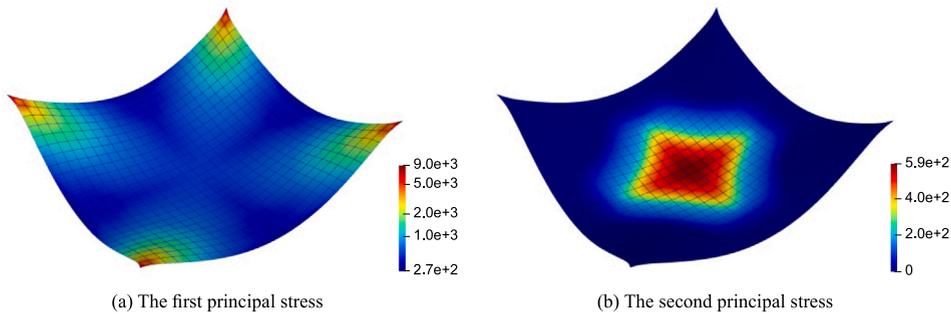

(a) The first principal stress                    (b) The second principal stress

**Fig. 22.** Distribution of the first and second principal stresses $\sigma_{1/2}$ (Pa) with $\nu = 0.3$.

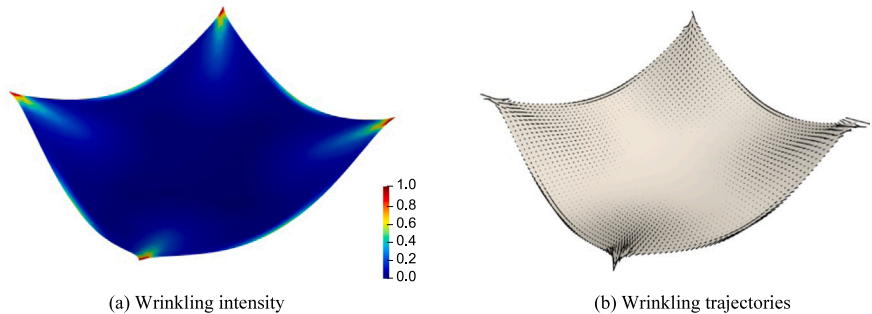

(a) Wrinkling intensity                    (b) Wrinkling trajectories

**Fig. 23.** Illustration of the wrinkled zone in a hanging blanket with $\nu = 0.3$.

Table 3 demonstrates very good agreement between the newly proposed model and the reference in predicting the displacements at points $A$, $B$, and $M$, denoted as $u_x^A$, $u_x^B$ and $u_z^M$, respectively, as well as the first principal stress $\sigma_1^M$ at the middle point $M$ when considering a Poisson's ratio $\nu = 0$. The results in Table 4 exhibit distinguishable variations between the two models for a nonzero Poisson's ratio. This discrepancy arises since the new wrinkling model is strain-based, while the wrinkling criterion used in the reference is mixed. However, the discrepancies remain within an acceptable range, indicating that the proposed model still provides reasonable predictions even for cases where the uniaxial tension condition is not strictly met.

**Table 3**
Comparison of the displacements and first principal stress with $\nu = 0$.

| Mesh ($25 \times 25$) | $u_z^M$ (m) | $u_x^A$ (m) | $u_x^B$ (m) | $\sigma_1^M$ (MPa) |
|---|---|---|---|---|
| Reference | −0.28949 | −0.03661 | −0.01830 | 637.68 |
| New model | −0.28949 | −0.03661 | −0.01830 | 611.68 |

**Table 4**
Comparison of the displacements and first principal stress with $\nu = 0.3$.

| Mesh ($25 \times 25$) | $u_z^M$ (m) | $u_x^A$ (m) | $u_x^B$ (m) | $\sigma_1^M$ (MPa) |
|---|---|---|---|---|
| Reference | −0.28328 | −0.03406 | −0.01703 | 642.66 |
| New model | −0.29531 | −0.03278 | −0.01639 | 586.99 |

The convergence performances of the two models are compared across various Poisson's ratios in Fig. 24. The results indicate that the new model exhibits significantly superior convergence compared to the reference. Specifically, the new wrinkling model achieves convergence within approximately twenty iterations, whereas the reference requires over one hundred iterations, particularly when the Poisson's ratio is zero. That attributes the advantages of the new model, which is consistently derived from the strain energy density. These findings underscore the improved convergence behavior of the proposed model, and the substantial reduction in the number of iterations needed for convergence demonstrates the advantages and practicality of the proposed model in simulating and analyzing wrinkling behavior in membrane structures.

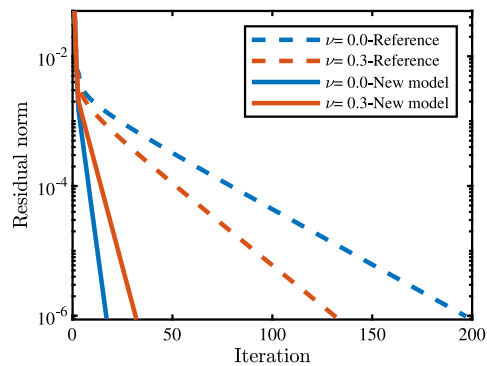

**Fig. 24.** Comparison of the convergence behavior of two models during the last step with various Poisson's ratios.

## 6. Conclusion

We have presented a variationally consistent membrane wrinkling model based on the spectral decomposition of the strain tensor, which can perform well in wrinkled membrane analysis and improve convergence issues caused by the wrinkling phenomenon. In this wrinkling model, we decompose the strain tensor into positive and negative components based on the eigenvalues. It allows us to split the strain energy density into positive and negative parts. According to the tension-field theory, we only consider the positive strain energy density to obtain zero compressive stiffness. To prevent self-penetration and singular stiffness matrices, a residual compressive stiffness can be assigned by multiplying the negative strain energy density by a very small factor, which is chosen empirically. We consistently derive the modified strain energy density with respect to strain variables to determine new stress tensor and constitutive tensor formulations. In addition, the eigenvector directions are considered to track the wrinkling directions.

To assess the effectiveness of the proposed method, we conducted a series of tests on analytical, numerical and experimental benchmarks for membrane wrinkling problems. We have employed this model in standard finite element and isogeometric analysis formulations, showing its generality. The results of these tests indicate that the newly proposed wrinkling model performs well in accurately predicting the mechanical responses of wrinkled thin membranes and shows a very good convergence behavior.

One possible direction for future research on the proposed approach is its extension to hyperelastic material models since the current model presented in this study is limited to isotropic linear elasticity. Another future research direction is the simultaneous spectral decomposition of the strain and stress tensors into positive and negative components. Thus, the wrinkling criterion would rely no longer solely on strain but on a widely used mixed stress–strain criterion.

**CRediT authorship contribution statement**

**Daobo Zhang:** Writing – original draft, Visualization, Validation, Software, Methodology, Investigation, Conceptualization. **Josef Kiendl:** Writing – review & editing, Visualization, Validation, Supervision, Software, Methodology, Investigation, Conceptualization.

**Declaration of competing interest**

The authors declare that they have no known competing financial interests or personal relationships that could have appeared to influence the work reported in this paper.

**Data availability**

Data will be made available on request.

**Appendix A**

In the following, the derivation of the positive and negative parts of the constitutive tensor are presented, following Nilsen [61]. For convenience, Eqs. (15) and (17) are repeated here

$$\mathbf{S}^{\pm} = \left( \lambda - \frac{\lambda^2}{\lambda + 2\mu} \right) \langle \operatorname{tr}(\mathbf{E}) \rangle^{\pm} \mathbf{I} + 2\mu \mathbf{E}^{\pm}, \tag{A.1}$$

$$\mathbb{C}^{\pm} = \frac{\partial \mathbf{S}^{\pm}(\mathbf{E})}{\partial \mathbf{E}}. \tag{A.2}$$

For a better clarity, the two terms in Eq. (A.1) are derived separately with respect to $\mathbf{E}$ in the following. For the first term, we get

$$\frac{\partial}{\partial \mathbf{E}} \left( \langle \operatorname{tr}(\mathbf{E}) \rangle^{\pm} \mathbf{I} \right) = H \left( \pm \operatorname{tr}(\mathbf{E}) \right) \frac{\partial}{\partial \mathbf{E}} \left( \operatorname{tr}(\mathbf{E}) \mathbf{I} \right), \tag{A.3}$$

where $H(x)$ is the Heaviside function. A fourth-order tensor $\mathbb{J}$ is introduced to define the derivative of the product of the trace of the strain tensor $\operatorname{tr}(\mathbf{E})$ and the identity tensor $\mathbf{I}$ as:

$$\left( \frac{\partial}{\partial \mathbf{E}} \left( \operatorname{tr}(\mathbf{E}) \mathbf{I} \right) \right)_{\alpha\beta\gamma\delta} := \mathbb{J} = \frac{\partial \left( \operatorname{tr}(\mathbf{E}) \mathbf{I} \right)_{\alpha\beta}}{\partial E_{\gamma\delta}} = \begin{cases} 1 & \text{if } \alpha = \beta \text{ and } \gamma = \delta \\ 0 & \text{else} . \end{cases} \tag{A.4}$$

Thereby, the Eq. (A.3) can be simply rewritten as:

$$\frac{\partial}{\partial \mathbf{E}} \left( \langle \operatorname{tr}(\mathbf{E}) \rangle^{\pm} \mathbf{I} \right) = H \left( \pm \operatorname{tr}(\mathbf{E}) \right) \mathbb{J} . \tag{A.5}$$

Next, we address the derivative of the second term in Eq. (15), which involves the so-called projection operators $\mathbb{P}_E^{\pm}$ [62]. These operators are fourth-order tensors and are defined as:

$$\mathbb{P}_E^{\pm} = \frac{\partial \mathbf{E}^{\pm}}{\partial \mathbf{E}} = \frac{\partial}{\partial \mathbf{E}} \left( \langle E_\alpha \rangle^{\pm} \mathbf{n}_\alpha \otimes \mathbf{n}_\alpha \right) . \tag{A.6}$$

The second-order tensors $\mathbf{M}_\alpha := \mathbf{n}_\alpha \otimes \mathbf{n}_\alpha$ are defined as the eigenvalue bases related to the eigenvalue $E_\alpha$. In the case of distinct eigenvalues $E_\alpha \neq E_\beta$, the classical results can be found in [63] and expressed as:

$$\frac{\partial E_\alpha}{\partial \mathbf{E}} = \mathbf{M}_\alpha, \qquad \frac{\partial \mathbf{M}_\alpha}{\partial \mathbf{E}} = \sum_{\alpha \neq \beta}^{2} \frac{1}{2 \left( E_\alpha - E_\beta \right)} \left( \mathbb{G}_{\alpha\beta} + \mathbb{G}_{\beta\alpha} \right), \tag{A.7}$$

where the fourth-order tensor operator $\mathbb{G}_{\alpha\beta} = \hat{\mathbb{G}}(\mathbf{M}_\alpha, \mathbf{M}_\beta)$ with the coordinates representation is introduced as:

$$\left( \mathbb{G}_{\alpha\beta} \right)_{\gamma\delta\epsilon\zeta} := (\mathbf{M}_\alpha)_{\gamma\epsilon} (\mathbf{M}_\beta)_{\delta\zeta} + (\mathbf{M}_\alpha)_{\gamma\zeta} (\mathbf{M}_\beta)_{\delta\epsilon}. \tag{A.8}$$

With this notation we can explicitly express the projection tensors $\mathbb{P}_E^{\pm} = \partial \mathbf{E}^{\pm}/\partial \mathbf{E}$. The components of the projection tensors are given by:

$$\begin{aligned} \left( \frac{\partial \mathbf{E}^{\pm}}{\partial \mathbf{E}} \right)_{\gamma\delta\epsilon\zeta} &= \frac{\partial \left( \mathbf{E}^{\pm} \right)_{\gamma\delta}}{\partial E_{\epsilon\zeta}} \\ &= \frac{\partial}{\partial E_{\epsilon\zeta}} \left( \langle E_\alpha \rangle^{\pm} \left( \mathbf{M}_\alpha \right)_{\gamma\delta} \right) \\ &= \frac{\partial \langle E_\alpha \rangle^{\pm}}{\partial E_{\epsilon\zeta}} \left( \mathbf{M}_\alpha \right)_{\gamma\delta} + \langle E_\alpha \rangle^{\pm} \frac{\partial \left( \mathbf{M}_\alpha \right)_{\gamma\delta}}{\partial E_{\epsilon\zeta}} \\ &= H \left( \pm E_\alpha \right) \left( \mathbf{M}_\alpha \right)_{\epsilon\zeta} \left( \mathbf{M}_\alpha \right)_{\gamma\delta} + \sum_{\alpha \neq \beta}^{2} \frac{\langle E_\alpha \rangle^{\pm}}{2 \left( E_\alpha - E_\beta \right)} \left( \left( \mathbb{G}_{\alpha\beta} \right)_{\gamma\delta\epsilon\zeta} + \left( \mathbb{G}_{\beta\alpha} \right)_{\gamma\delta\epsilon\zeta} \right). \end{aligned} \tag{A.9}$$

Finally, we can express the positive and negative parts of the constitutive tensor as

$$\mathbb{C}^{\pm} = \left( \lambda - \frac{\lambda^2}{\lambda + 2\mu} \right) H \left( \pm \operatorname{tr}(\mathbf{E}) \right) \mathbb{J} + 2\mu \left( H \left( \pm E_\alpha \right) \mathbb{Q}_\alpha + \sum_{\alpha \neq \beta}^{2} \frac{\langle E_\alpha \rangle^{\pm}}{2 \left( E_\alpha - E_\beta \right)} \left( \mathbb{G}_{\alpha\beta} + \mathbb{G}_{\beta\alpha} \right) \right). \tag{A.10}$$

## Appendix B

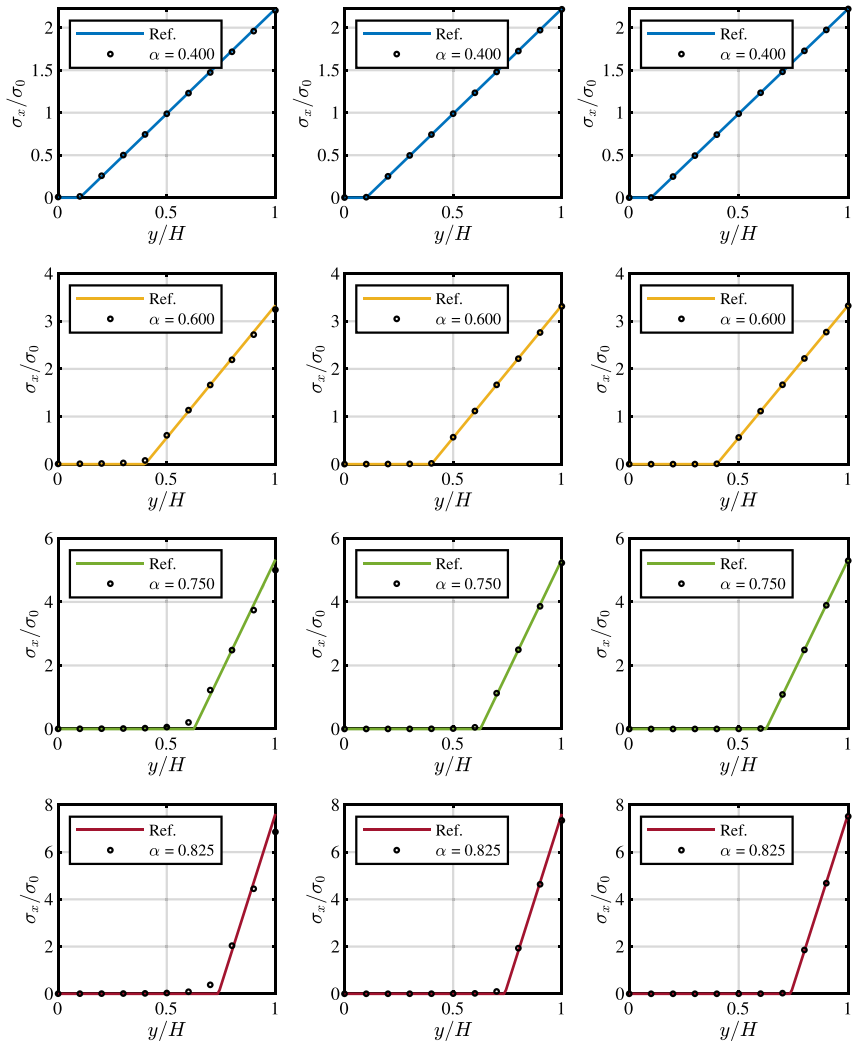

**Fig. B.25.** Normal stress $\sigma_x$ distribution along the height $H$ computed by the proposed model with Poisson's ratio $\nu = 0$, polynomial degree $p = 1$, and $\alpha = 2M/PH$. In the left, middle, and right columns, mesh refinement is carried out as $11 \times 5$, $23 \times 11$, and $43 \times 21$, respectively.

## References

[1] Y. Miyazaki, Wrinkle/slack model and finite element dynamics of membrane, Internat. J. Numer. Methods Engrg. 66 (7) (2006) 1179–1209.

[2] X. Liu, C.H. Jenkins, W.W. Schur, Large deflection analysis of pneumatic envelopes using a penalty parameter modified material model, Finite Elem. Anal. Des. 37 (3) (2001) 233–251.

[3] R. Rossi, M. Lazzari, R. Vitaliani, E. Oñate, Simulation of light-weight membrane structures by wrinkling model, Internat. J. Numer. Methods Engrg. 62 (15) (2005) 2127–2153.

[4] X.W. Feng, Q.S. Yang, S.-s. Law, Wrinkled membrane element based on the wrinkling potential, Int. J. Solids Struct. 51 (21–22) (2014) 3532–3548.

[5] M. Stein, J.M. Hedgepeth, Analysis of Partly Wrinkled Membranes, National Aeronautics and Space Administration, 1961.

[6] P. Contri, B. Schrefler, A geometrically nonlinear finite element analysis of wrinkled membrane surfaces by a no-compression material model, Commun. Appl. Numer. Methods 4 (1) (1988) 5–15.

[7] R.K. Miller, J.M. Hedgepeth, An algorithm for finite element analysis of partly wrinkled membranes, AIAA J. 20 (12) (1982) 1761–1763.

[8] R.K. Miller, J.M. Hedgepeth, V.I. Weingarten, P. Das, S. Kahyai, Finite element analysis of partly wrinkled membranes, in: Advances and Trends in Structures and Dynamics, Elsevier, 1985, pp. 631–639.

[9] H. Wagner, Flat Sheet Metal Girders with Very Thin Metal Web. Part I: General Theories and Assumptions, Technical report, 1931.

[10] E. Reissner, On tension field theory, in: Proc. of the 5th Int. Congr. for Applied Mechanics Harvard Univ. & MIT, 1938, pp. 88–92.

[11] E. Mansfield, Tension field theory, a new approach which shows its duality with inextensional theory, in: Applied Mechanics: Proceedings of the Twelfth International Congress of Applied Mechanics, Stanford University, 1968, pp. 26–31, pages 305–320. Springer, 1969.

[12] C.H. Wu, Nonlinear wrinkling of nonlinear membranes of retolution, J. Appl. Mech. Trans. ASME 45 (3) (1978) 533–538.

[13] E. Cerda, L. Mahadevan, Geometry and physics of wrinkling, Phys. Rev. Lett. 90 (7) (2003) 074302.

[14] W. Wong, S. Pellegrino, Wrinkled membranes II: analytical models, J. Mech. Mater. Struct. 1 (1) (2006) 27–61.

[15] E. Puntel, L. Deseri, E. Fried, Wrinkling of a stretched thin sheet, J. Elasticity 105 (2011) 137–170.

[16] W. Wong, S. Pellegrino, Wrinkled membranes III: numerical simulations, J. Mech. Mater. Struct. 1 (1) (2006) 63–95.

[17] M. Taylor, K. Bertoldi, D.J. Steigmann, Spatial resolution of wrinkle patterns in thin elastic sheets at finite strain, J. Mech. Phys. Solids 62 (2014) 163–180.

[18] A. Jarasjarungkiat, R. Wüchner, K.-U. Bletzinger, A wrinkling model based on material modification for isotropic and orthotropic membranes, Comput. Methods Appl. Mech. Engrg. 197 (6–8) (2008) 773–788.

[19] T. Iwasa, M. Natori, K. Higuchi, Evaluation of tension field theory for wrinkling analysis with respect to the post-buckling study, J. Appl. Mech. 71 (4) (2004) 532–540.

[20] H.M. Verhelst, M. Möller, J. Den Besten, A. Mantzaflaris, M.L. Kaminski, Stretch-based hyperelastic material formulations for isogeometric Kirchhoff–Love shells with application to wrinkling, Comput. Aided Des. 139 (2021) 103075.

[21] F.G. Flores, E. Oñate, Wrinkling and folding analysis of elastic membranes using an enhanced rotation-free thin shell triangular element, Finite Elem. Anal. Des. 47 (9) (2011) 982–990.

[22] C. Fu, H-H. Dai, F. Xu, Computing wrinkling and restabilization of stretched sheets based on a consistent finite-strain plate theory, Comput. Methods Appl. Mech. Engrg. 384 (2021) 113986.

[23] D. Steigmann, Tension-field theory, Proc. R. Soc. A 429 (1876) (1990) 141–173.

[24] D.G. Roddeman, J. Drukker, C.W.J. Oomens, J.D. Janssen, The wrinkling of thin membranes: Part I—theory, J. Appl. Mech. 54 (4) (1987) 884–887.

[25] D.G. Roddeman, J. Drukker, C.W.J. Oomens, J.D. Janssen, The wrinkling of thin membranes: Part II—numerical analysis, J. Appl. Mech. 54 (4) (1987) 888–892.

[26] C.H. Wu, T.R. Canfield, Wrinkling in finite plane-stress theory, Quart. Appl. Math. 39 (2) (1981) 179–199.

[27] K. Lu, M. Accorsi, J. Leonard, Finite element analysis of membrane wrinkling, Int. J. Numer. Methods Eng. 50 (5) (2001) 1017–1038.

[28] H. Schoop, L. Taenzer, J. Hornig, Wrinkling of nonlinear membranes, Comput. Mech. 29 (2002) 68–74.

[29] K. Nakashima, M.C. Natori, Efficient modification scheme of stress–strain tensor for wrinkled membranes, AIAA J. 43 (1) (2005) 206–215.

[30] K. Nakashino, M. Natori, Three-dimensional analysis of wrinkled membranes using modification scheme of stress–strain tensor, AIAA J. 44 (7) (2006) 1498–1504.

[31] J. Hornig, H. Schoop, Closed form analysis of wrinkled membranes with linear stress–strain relation, Comput. Mech. 30 (4) (2003) 259–264.

[32] K. Nakashino, A. Nordmark, A. Eriksson, Geometrically nonlinear isogeometric analysis of a partly wrinkled membrane structure, Comput. Struct. 239 (2020) 106302.

[33] S. Kang, S. Im, Finite element analysis of wrinkling membranes, J. Appl. Mech. 64 (2) (1997) 263–269.

[34] S. Kang, S. Im, Finite element analysis of dynamic response of wrinkling membranes, Comput. Methods Appl. Mech. Engrg. 173 (1–2) (1999) 227–240.

[35] T. Raible, K. Tegeler, S. Löhnert, P. Wriggers, Development of a wrinkling algorithm for orthotropic membrane materials, Comput. Methods Appl. Mech. Engrg. 194 (21–24) (2005) 2550–2568.

[36] K. Woo, H. Igawa, C. Jenkins, Analysis of wrinkling behavior of anisotropic membrane, Comput. Model. Eng. Sci. 6 (2004) 397–408.

[37] T. Akita, K. Nakashino, M. Natori, K. Park, A simple computer implementation of membrane wrinkle behaviour via a projection technique, Int. J. Numer. Methods Eng. 71 (10) (2007) 1231–1259.

[38] A. Jarasjarungkiat, R. Wüchner, K.-U. Bletzinger, Efficient sub-grid scale modeling of membrane wrinkling by a projection method, Comput. Methods Appl. Mech. Engrg. 198 (9–12) (2009) 1097–1116.

[39] H. Le Meitour, G. Rio, H. Laurent, A. Lectez, P. Guigue, Analysis of wrinkled membrane structures using a plane stress projection procedure and the dynamic relaxation method, Int. J. Solids Struct. 208 (2021) 194–213.

[40] A.C. Pipkin, The relaxed energy density for isotropic elastic membranes, IMA J. Appl. Math. 36 (1) (1986) 85–99.

[41] E. Haseganu, D. Steigmann, Analysis of partly wrinkled membranes by the method of dynamic relaxation, Comput. Mech. 14 (6) (1994) 596–614.

[42] M. Epstein, M.A. Forcinito, Anisotropic membrane wrinkling: theory and analysis, Int. J. Solids Struct. 38 (30–31) (2001) 5253–5272.

[43] J. Mosler, A novel variational algorithmic formulation for wrinkling at finite strains based on energy minimization: application to mesh adaption, Comput. Methods Appl. Mech. Engrg. 197 (9–12) (2008) 1131–1146.

[44] J. Mosler, F. Cirak, A variational formulation for finite deformation wrinkling analysis of inelastic membranes, Comput. Methods Appl. Mech. Engrg. 198 (27–29) (2009) 2087–2098.

[45] C. Miehe, M. Hofacker, F. Welschinger, A phase field model for rate-independent crack propagation: Robust algorithmic implementation based on operator splits, Comput. Methods Appl. Mech. Engrg. 199 (45–48) (2010) 2765–2778.

[46] J. Kiendl, M. Ambati, L. De Lorenzis, H. Gomez, A. Reali, Phase-field description of brittle fracture in plates and shells, Comput. Methods Appl. Mech. Engrg. 312 (2016) 374–394.

[47] J.A. Cottrell, T.J. Hughes, Y. Bazilevs, Isogeometric Analysis: Toward Integration of CAD and FEA, John Wiley & Sons, 2009.

[48] T.J. Hughes, J.A. Cottrell, Y. Bazilevs, Isogeometric analysis: cad, finite elements, nurbs, exact geometry and mesh refinement, Comput. Appl. Mech. Engrg. 194 (39–41) (2005) 4135–4195.

[49] J. Kiendl, M.-C. Hsu, M.C. Wu, A. Reali, Isogeometric Kirchhoff–Love shell formulations for general hyperelastic materials, Comput. Methods Appl. Mech. Engrg. 291 (2015) 280–303.

[50] C. Liu, Q. Tian, D. Yan, H. Hu, Dynamic analysis of membrane systems undergoing overall motions, large deformations and wrinkles via thin shell elements of ancf, Comput. Methods Appl. Mech. Engrg. 258 (2013) 81–95.

[51] X.W. Feng, J. Ma, S.-s. Law, Q. Yang, Numerical analysis of wrinkle-influencing factors of thin membranes, Int. J. Solids Struct. 97 (2016) 458–474.

[52] H. Ding, B. Yang, M. Lou, H. Fang, New numerical method for two-dimensional partially wrinkled membranes, AIAA J. 41 (1) (2003) 125–132.

[53] D.G. Jeong, B.M. Kwak, Complementarity problem formulation for the wrinkled membrane and numerical implementation, Finite Elem. Anal. Des. 12 (2) (1992) 91–104.

[54] W. Wong, S. Pellegrino, Wrinkled membranes I: experiments, J. Mech. Mater. Struct. 1 (1) (2006) 3–25.

[55] A. Diaby, C. Wielgosz, et al., Buckling and wrinkling of prestressed membranes, Finite Elem. Anal. Des. 42 (11) (2006) 992–1001.

[56] E.-S. Lee, S.-K. Youn, Finite element analysis of wrinkling membrane structures with large deformations, Finite Elem. Anal. Des. 42 (8–9) (2006) 780–791.

[57] H. Le Meitour, G. Rio, H. Laurent, A. Lectez, P. Guigue, Analysis of wrinkled membrane structures using a plane stress projection procedure and the dynamic relaxation method, Int. J. Solids Struct. 208 (2021) 194–213.

[58] A.J. Gil, J. Bonet, Finite element analysis of partly wrinkled reinforced prestressed membranes, Comput. Mech. 40 (2007) 595–615.

[59] J. Lu, C. Zheng, Dynamic cloth simulation by isogeometric analysis, Comput. Methods Appl. Mech. Engrg. 268 (2014) 475–493.

[60] A.J. Herrema, E.L. Johnson, D. Proserpio, M.C. Wu, J. Kiendl, M.-C. Hsu, Penalty coupling of non-matching isogeometric Kirchhoff–Love shell patches with application to composite wind turbine blades, Comput. Methods Appl. Mech. Engrg. 346 (2019) 810–840.

[61] O.K. Nilsen, Simulation of Crack Propagation using Isogeometric Analysis Applied with Nurbs and Lr B-Splines Master's thesis, Institutt for matematiske fag, 2012.

[62] V. Lubarda, D. Krajcinovic, S. Mastilovic, Damage model for brittle elastic solids with unequal tensile and compressive strengths, Eng. Fract. Mech. 49 (5) (1994) 681–697.

[63] C. Miehe, M. Lambrecht, Algorithms for computation of stresses and elasticity moduli in terms of seth–hill's family of generalized strain tensors, Commun. Numer. Methods Eng. 17 (5) (2001) 337–353.